\documentclass[11pt, a4paper]{article}
\usepackage[utf8]{inputenc}
\usepackage{mlmodern}
\usepackage[T1]{fontenc}
\usepackage[english]{babel}
\usepackage[hidelinks,hypertexnames=false]{hyperref}
\usepackage{amssymb}
\usepackage{mathtools}
\usepackage{mathrsfs}
\usepackage{amsthm}
\usepackage{orcidlink}
\usepackage{graphicx}
\usepackage{enumerate}
\usepackage{microtype}
\usepackage{thmtools}
\usepackage[capitalise]{cleveref}
\usepackage[backend=biber,style=alphabetic,sorting=ynt]{biblatex}
\usepackage{comment}

\usepackage[left]{lineno}
\AtBeginDocument{}

\DeclareSourcemap{
  \maps[datatype=bibtex]{
    \map{
      \step[fieldset=howpublished, null]
      \step[fieldset=isbn,null]
      \step[fieldset=issn,null]
      \step[fieldset=language,null]
    }
  }
}

\makeatletter
\def\th@plain{%
  \thm@notefont{}%
  \itshape%
}
\def\th@definition{%
  \thm@notefont{}%
  \normalfont%
}
\makeatother

\numberwithin{equation}{section}

\theoremstyle{definition}
\newtheorem{definition}{Definition}[section]
\crefname{definition}{Definition}{Definitions}

\newtheorem{remark}[definition]{Remark}
\crefname{remark}{Remark}{Remarks}

\newtheorem{lemma}[definition]{Lemma}
\crefname{lemma}{Lemma}{Lemmata}

\newtheorem{proposition}[definition]{Proposition}
\crefname{proposition}{Proposition}{Propositions}

\newtheorem{theorem}[definition]{Theorem}
\crefname{theorem}{Theorem}{Theorems}

\crefname{corollary}{Corollary}{Corollaries}

\crefname{conjecture}{Conjecture}{Conjectures}

\crefname{counterexample}{Counterexample}{Counterexamples}

\DeclareMathOperator{\re}{Re}
\DeclareMathOperator{\im}{Im}
\DeclareMathOperator{\ran}{ran}
\DeclareMathOperator{\dom}{dom}
\DeclareMathOperator{\supp}{supp}
\DeclareMathOperator{\dist}{dist}

\newcommand{\diff}{\, \mathrm{d}}
\newcommand{\id}{\mathrm{I}}

\newcommand{\NN}{\mathbb{N}}
\newcommand{\CC}{\mathbb{C}}
\newcommand{\RR}{\mathbb{R}}

\makeatletter
\newcommand{\smallbullet}{}%
\DeclareRobustCommand\smallbullet{%
  \, \mathord{\mathpalette\smallbullet@{0.6}} \,
}
\newcommand{\smallbullet@}[2]{%
  \vcenter{\hbox{\scalebox{#2}{$\m@th#1\bullet$}}}%
}
\makeatother

\let\oldforall\forall%
\let\forall\undefined%
\DeclareMathOperator{\forall}{\oldforall}

\title{On the sufficiency of the Hautus test for exact observability
  of analytic normal semigroups\thanks{This work was funded by the Deutsche
    Forschungsgemeinschaft
(DFG, German Research Foundation) – Project-ID 531152215 – CRC 1701.}}
\author{Jannik Daun\orcidlink{0009-0006-9455-972X}\thanks{University
    of Wuppertal, School of
    Mathematics and Natural Sciences, Gaußstraße 20, 42119
    Wuppertal, Germany,
  (jdaun@uni-wuppertal.de, bjacob@uni-wuppertal.de).}
\and Birgit Jacob\orcidlink{0000-0002-8775-4424}\footnotemark[2]}
\date{\today}
\begin{document}
\maketitle
\begin{abstract}
  We prove that the infinite-dimensional Hautus test implies
  exact observability
  for normal analytic exponentially stable strongly continuous
  semigroups and admissible
  observation operators.

  We show how to verify the Hautus test to conclude exact
  observability for a diffusion
  equation with a singular potential and derivative observation as well as a
  multiplication operator with high-pass frequency domain measurements.
\end{abstract}
\renewcommand{\abstractname}{Key words}
\begin{abstract}
  Exact observability, Hautus
  test, Russell--Weiss
  conjecture, normal generator, normal
  semigroup, analytic semigroup, sectorial operator, admissible
  observation operators
\end{abstract}
\renewcommand{\abstractname}{MSC codes}
\begin{abstract}
  93B07, 93C25, 93C05, 47D06, 93B28
\end{abstract}
\section{Introduction and main results}
We consider a system $(A,C)$, where $A$ is the generator of an
exponentially stable strongly continuous semigroup $T$ on the Hilbert
space $X$ and $C : X_1 \to Y$ is an admissible observation operator.
Here $Y$ is another Hilbert space and $X_1$ is $\dom (A)$ with the
norm $\| v\|_1 \coloneq \| A v \|$.

In applications $A$ is often a linear partial differential operator
on an appropriate $L^2$-space,
for example $X= L^2(\RR^n)$ and $A =  \Delta - \id$.
For $v \in X_1$ the function $x : [0,\infty) \to X_1$, $x(t) \coloneq
T(t) v$ is the unique $x \in C^1([0,\infty),X)$ which satisfies the
initial value problem
\begin{equation*}
  \begin{split}
    \dot x(t) &= A x(t), \ t \in [0,\infty), \\
    x(0) &= v.
  \end{split}
\end{equation*}
The function $x$ encodes the state of the system at different times.
The system's output $y \in L^2([0,\infty),Y)$ is defined by $y (t)
\coloneq C x(t)$.
We want to investigate if the initial state $v$ can be
reconstructed continuously in the $X$-norm from the measured output $y$.
This property is encoded in the following definition.
The system $(A,C)$ is called \emph{exactly observable} if and only if
\begin{equation*}
  0 < \inf_{ v \in X_1, \ \| v\|=1 } \int_0^\infty \| C T(t) v
  \|^2 \diff \lambda (t).
\end{equation*}
Showing exact observability directly from the definition can be
challenging, especially when explicit expressions for the semigroup
$T$ are not available.
For this reason it is of great interest to find an equivalent
characterization of exact observability
in terms of the system's generator $A$ and the observation operator $C$.
A candidate for such
an equivalent characterization is the Hautus test, defined as follows:
The system $(A,C)$ is said to satisfy the (infinite-dimensional)
\emph{Hautus test} if there exists $m \in (0,\infty)$ such that for
all $z \in \CC_-$ and $v \in X_1$ we have
\begin{equation*}
  m \| v\|^2    \leq | \re z |^{-2} \| (z \id -A) v \|^2 + |\re z
  |^{-1} \| C v\|^2.
\end{equation*}
Russell and Weiss~\cite{RussellWeissGeneralCondExactObs} introduced
the Hautus test, showed that it is necessary for exact observability
and conjectured that it is also sufficient.
The Hautus test is known to be sufficient for exact observability, if
$A$ is bounded~\cite{RussellWeissGeneralCondExactObs}, if $A = B -\id$
with $B$ skew-adjoint~\cite{miller_controllability_2005}, if $T$ is
diagonal and $Y$
finite-dimensional~\cite{jacobZwartExactObsFinDimOutput}, if $T$
extends to a normal group~\cite{jacobZwartHautusForExpStable}, and if
$A$ is self-adjoint~\cite{daun2026hautustestexactobservability}.
In general the Hautus test does not imply exact observability.
A counterexample was given
in~\cite{jacob_counterexamples_2004}. In this counterexample $Y =
\CC$ and $T$ is analytic, but not similar to a contraction semigroup.
Subsequently it was conjectured that the Hautus test is sufficient for exact
observability if $T$ is similar to
a contraction semigroup~\cite{jacob_zwart_blondel_unsolved_2004}.
Recently we constructed a counterexample to this conjecture in which
$A$ is normal and
has an orthonormal basis of
eigenvectors~\cite{daun2026hautustestexactobservability}. The
normality of $A$ and the exponential stability of $T$ imply that $T$
is a contraction semigroup.

In summary normality of $A$ and analyticity of $T$ are
individually insufficient for the Hautus test to imply exact
observability, while self-adjointness of $A$ is sufficient.
In~\cite{daun2026hautustestexactobservability} we asked if the
combination of analyticity and normality is also sufficient.
The affirmative answer to this question is the main new result in this article.
\begin{theorem}\label{introduction:mainThm}
  Let $A$ be the generator of an
  exponentially stable strongly continuous semigroup $T$ on the Hilbert
  space $X$ and $C : X_1 \to Y$ an admissible observation operator
  taking values in the Hilbert space $Y$.
  Assume that $T$ is analytic and normal. Then the Hautus test
  implies exact observability.
\end{theorem}
We now briefly outline the proof of \cref{introduction:mainThm}.
Let $E$ be the spectral measure of $A$.
We fix a well chosen bounded annulus $\Omega \subset \CC$. For $t \in
(0,\infty)$ we let
$P_t \coloneq E (t^{-1}\Omega)$, $V_t \coloneq \ran P_t$, $A_t
\coloneq t A|_{V_t}$, and $C_t \coloneq
\sqrt{t} C|_{V_t}$.
We consider the Hilbert space direct sum of the systems $(A_t,C_t)$
for $t\in (0,\infty)$.
Using a scaling symmetry in the Hautus test we show that the direct
sum system also satisfies the Hautus test.
Moreover we show that its generator is bounded and that the direct
sum system is exponentially stable.
Applying that the Hautus test is sufficient for exact observability
in the bounded generator case to the direct sum system yields a
uniform observability estimate for all systems.
To be precise we show the existence of $m_o \in (0,\infty)$ and $L
\in (0,\infty)$ such that
for all $t \in
(0,\infty)$ and $v \in V_t$
it is true that
\begin{equation}
  \label{eq:intro0}
  m_o \|v\|^2 \leq t \int_0^L \| C T(t s) v\|^2 \diff \lambda (s).
\end{equation}
The remaining difficulty is to glue the local observability
estimate on each space $V_t$ into a global estimate.

To this end we consider the function $f : \CC
\to \CC$, $f(z) \coloneq {( e^z (1-e^z))}^N$ for some sufficiently large
$N\in \NN$.
We show that there are
$\kappa, c,d \in (0,\infty)$ with $m_o c
- 2 d>0$ and such that
for all $v \in X_1$ we have
\begin{equation}
  \label{eq:intro1}
  c \| v\|^2 \leq \int_0^\infty t^{-1}\| f (t A) P_t v\|^2 \diff \lambda (t),
\end{equation}
\begin{equation}
  \label{eq:intro2}
  I_{1,v} \coloneq \int_0^\infty \int_0^L  \|C T(t s)
  f(t A) v\|^2
  \diff \lambda (s) \diff
  \lambda (t)\leq  \kappa \int_0^\infty  \| C
  T(t) v\|^2 \diff
  \lambda (t),
\end{equation}
and
\begin{equation}
  \label{eq:intro3}
  I_{2,v} \coloneq  \int_0^\infty \int_0^L
  \| C T(t s) f(t A)
  (I- P_t) v\|^2
  \diff \lambda (s) \diff
  \lambda (t) \leq d \| v\|^2.
\end{equation}
Let $v \in X_1$.
From~\eqref{eq:intro1},~\eqref{eq:intro0} and the multiplication
theorem for the functional calculus we get
\begin{equation}
  \begin{split}
    m_o c \|v\|^2 &\leq m_o \int_0^\infty t^{-1}\| f(t A) P_t
    v\|^2 \diff
    \lambda (t)  \\
    &\leq   \int_0^\infty \int_0^L  \|C T(ts) f(t A) P_t v\|^2
    \diff \lambda (s) \diff
    \lambda (t)  \\
    &\leq 2 ( I_{1,v} + I_{2,v}).
  \end{split}
\end{equation}
Applying~\eqref{eq:intro2} and~\eqref{eq:intro3} yields
\begin{equation*}
  \underbrace{(m_o c - 2  d)}_{>0} \|v\|^2 \leq 2 \kappa
  \int_0^\infty  \| C T(t)
  v\|^2 \diff \lambda (t),
\end{equation*}
which is the exact observability of $(A,C)$.

To illustrate the utility of \cref{introduction:mainThm} we give two example
applications.
As our first example we apply \cref{introduction:mainThm} to a normal
and generally not self-adjoint
multiplication semigroup on $L^2((0,\infty))$ whose generator
has continuous spectrum. In this example our observation is given by
high-pass frequency measurements of the weighted state.

In the second example we consider a diffusion equation with a singular potential
consisting of disjointly supported spikes accumulating at the boundary of the
underlying domain.
The observation operator is given by the
derivative.
We verify the Hautus estimate without computing the spectrum or the
spectral measure of the generator.
In \cref{example2:generalizationRemark} we also highlight a general
mechanism for proving exact observability
of diffusion equations with singular potentials and gradient
observation via the Hautus test,
which is applicable beyond our specific example.

Both of our examples feature an unbounded generator,
infinite dimensional output space $Y$, and fail the lower energy
dissipation comparison criterion for exact observability in
\cref{examples:RealPartEOSufficiencyEstimate}.

This article is organized as follows.
Preliminary results on the direct sum of systems which are needed in
the proof of \cref{introduction:mainThm} are collected in
\cref{section:uniformEO}.
Then the proof of \cref{introduction:mainThm} is carried out in
\cref{sec:sector}.
In \cref{sec:criteria} we collect energy dissipation comparison criteria for
exact observability and admissibility and a sufficient condition for
the Hautus test that
only involves the real part of the
generator. These criteria are used when discussing the examples. The
multiplication operator
example is treated in \cref{sec:example1} and the diffusion equation
example in \cref{sec:example2}.
\subsection*{Notation}
We denote by $\lambda$ the Lebesgue measure on the Borel
$\sigma$-algebra of $\RR$.
By $\mathscr{B}(E)$ we denote the Banach algebra of bounded linear
operators $E \to E$ for some Banach space $E$.
We let $\CC_\pm \coloneq \{ z \in \CC \colon \pm \re z  >0 \}$, $\NN$ the
natural numbers including
$0$, and $\NN_+\coloneq \NN \setminus \{0\}$.
Throughout all vector spaces are assumed to be complex.

\section{Uniform exact observability}\label{section:uniformEO}
In this section we show that a family of systems with uniformly
bounded generators and uniform estimates for admissibility,
exponential stability and the Hautus test
is also uniformly exactly observable in a common finite time (see
\cref{uniformEO:mainThm}).

Let $J$ be an arbitrary non-empty set. For $j \in J$ let $X_j$ and
$Y_j$ be Hilbert spaces, let $A_j$ be a bounded linear operator on
$X_j$, and let $C_j:
X_j \to Y_j$ be a bounded linear operator.
For $j \in J$ let $T_j : [0,\infty) \to
\mathscr{B}(X_j)$, $T_j (t) \coloneq
e^{t A_j}$, i.e.~$T_j$ is the (strongly) continuous semigroup
generated by $A_j$.
Assume additionally that
\begin{enumerate}[(i)]
  \item $(A_j)_{j\in J}$ is uniformly bounded, i.e.~$\sup_{j \in J}
    \|A_j\| < \infty$.
  \item $(T_j)_{j\in J}$ is uniformly exponentially stable,
    i.e.~there exists $\kappa, \omega \in (0,\infty)$ such that for
    all $j \in J$ and $t \in [0,\infty)$ we have
    \begin{equation*}
      \| T_j(t)\| \leq \kappa e^{-\omega t}.
    \end{equation*}
  \item $(C_j)_{j\in J}$ is uniformly admissible, i.e.~there exists
    $M \in (0,\infty)$ such that for all $j \in J$ and $v \in X_j$ we have
    \begin{equation*}
      \int_0^\infty \|C_j T_j(t) v\|^2 \diff \lambda (t) \leq M \|v\|^2.
    \end{equation*}
\end{enumerate}
To define the direct sum $\oplus_{j \in J} C_j$ as a bounded operator
we need a uniform bound.
\begin{proposition}\label{uniformEO:CuniformBdd}
  We have
  \begin{equation*}
    \sup_{j \in J} \| C_j \| <\infty.
  \end{equation*}
\end{proposition}
\begin{proof}
  Let $b \coloneq \sup_{j \in J} \|A_j\| < \infty$.
  Let $j \in J$.
  For $h \in (0,\infty)$ define the bounded linear operator $Q_h$ on $X_j$ by
  \begin{equation*}
    Q_h \coloneq \frac{1}{h}\int_0^h T_j(s) \diff \lambda (s).
  \end{equation*}
  For every $h \in (0,\infty)$ we have that
  \begin{equation*}
    \begin{split}
      \| Q_h  - \id \| &\leq \frac{1}{h}\int_0^h \| T_j(s) -
      \id \|\diff \lambda (s) \\
      &\leq \frac{1}{h}\int_0^h \| A_j s \| \exp (\| A_j s\|) \diff
      \lambda (s) \\
      &\leq  b h  \exp ( b h).
    \end{split}
  \end{equation*}
  From now assume that $h$ is chosen such that $b h  \exp ( b h)< 1$.
  Since
  \begin{equation*}
    Q_h = Q_h -\id + \id = \id -(-Q_h+\id )
  \end{equation*}
  and $\|Q_h -\id \|<1$ we obtain from
  the geometric series that $Q_h $ is invertible and that
  $\|Q_h^{-1}\| \leq \frac{1}{1- bh \exp(bh)}$.
  Let $v \in X_j$.
  From Cauchy-Schwarz and the uniform admissibility of ${(C_j)}_{j \in
  J}$ we get
  \begin{equation*}
    \begin{split}
      \| C_j Q_h v\| &\leq \frac{1}{h} \int_0^h \| C_j T_j(s) v
      \| \diff \lambda (s) \\
      &\leq \frac{1}{\sqrt{h}} {\bigg(\int_0^h \| C_j T_j(s) v \|^2 \diff
      \lambda (s) \bigg)}^{1/2} \\
      &\leq \frac{1}{\sqrt{h}}  {\bigg(\int_0^\infty \| C_j T_j(s) v \|^2
          \diff \lambda (s)
      \bigg)}^{1/2} \\
      &\leq \frac{\sqrt{M}}{\sqrt{h}} \|v\| .
    \end{split}
  \end{equation*}
  Therefore for all $j \in J$ it is true that
  \begin{equation*}
    \| C_j \| \leq \| C_j Q_h \| \| Q_h^{-1}\|\leq
    \frac{\sqrt{M}}{\sqrt{h}} \frac{1}{1-bh \exp (bh)}. \qedhere
  \end{equation*}
\end{proof}
\begin{remark}
  The uniform admissibility of $(C_j)_{j \in J}$ and the property
  $\sup_{j \in J} \|C_j\|<\infty$ are actually equivalent by
  \cref{uniformEO:CuniformBdd} and the uniform exponential stability.
\end{remark}
\begin{theorem}
  \label{uniformEO:mainThm}
  Assume additionally that
  $(A_j)_{j\in J}$ and $(C_j)_{j\in J}$ satisfy the
  Hautus test uniformly, i.e.~there exists $m_h \in (0,\infty)$ such that for
  all $j \in J, v\in X_j$, and $z \in \CC_-$ it
  is true that
  \begin{equation*}
    m_h \|v\|^2 \leq \frac{1}{|\re z|} \| C_j v\|^2 + \frac{1}{{|\re
    z|}^2} \| (z\id - A_j) v\|^2.
  \end{equation*}
  Then $(A_j)_{j\in J}$ and $(C_j)_{j\in J}$ are uniformly exactly
  observable in finite time, i.e.~there exists $L \in (0,\infty)$,
  $m_o \in (0,\infty)$ such that for all $j \in J$ we have that for
  all $v \in X_j$ it is true that
  \begin{equation*}
    m_o \|v\|^2\leq \int_0^L \| C_j T_j(t) v\|^2 \diff \lambda (t).
  \end{equation*}
\end{theorem}
\begin{proof}
  Let $\mathcal{X}\coloneq \bigoplus_{j \in J} X_j$ (direct sum of
  Hilbert spaces) and for
  $j \in J$
  let $\iota_j : X_j \to \mathcal{X}$ be the embedding.
  Let $\mathcal{A}$ be the unique bounded linear operator on
  $\mathcal{X}$ such that for all $j \in J$ we have
  $ \mathcal{A} \iota_j =\iota_j  A_j$.
  The operator $\mathcal{A}$ exists, because $(A_j)_{j \in J}$ is
  uniformly bounded.
  Let $\mathcal{Y} \coloneq \bigoplus_{j \in J} Y_j$ and for
  $j \in J$
  let $\tilde{\iota}_j : Y_j \to \mathcal{Y}$ be the embedding.
  Let $\mathcal{C} :\mathcal{X} \to \mathcal{Y}$ be the unique bounded
  linear operator such that for all $j \in J$ we have
  $\mathcal{C}  \iota_j = \tilde{\iota}_j C_j$.
  The operator $\mathcal{C}$ exists, because $(C_j)_{j \in J}$ is
  uniformly bounded by \cref{uniformEO:CuniformBdd}.

  Let $\mathcal{T} : [0,\infty) \to \mathscr{B}(\mathcal{X})$, $\mathcal{T}
  (s)\coloneq \exp (s
  \mathcal{A} )$.
  Note that for all $s \in [0,\infty)$ and $j \in J$ we have
  $\mathcal{T}(s) \iota_j =
  \iota_j T_j(s)$.
  The uniform exponential stability of $(T_j)_{j \in J}$ implies that
  $\mathcal{T}$ is exponentially stable. Since $\mathcal{A}$ is
  bounded and $\mathcal{T}$ is exponentially stable we have that
  $\mathcal{C}$ is admissible (because it is bounded).
  By assumption $(\mathcal{A}, \mathcal{C})$
  satisfies the Hautus test.
  By~\cite[Theorem 1.7, p.~3]{RussellWeissGeneralCondExactObs} this implies that
  $(\mathcal{A}, \mathcal{C})$ is exactly observable.
  By the exponential stability infinite time exact observability implies
  finite time exact observability for some time~\cite[Proposition
  6.5.2, p.~194]{TucsnakWeiss}.
\end{proof}

\section{Sufficiency of the Hautus test for normal analytic
semigroups}\label{sec:sector}
In this section we prove \cref{introduction:mainThm}, i.e.~that the
Hautus test implies
exact observability for normal analytic semigroups and admissible
observation operators.
To do so we will use the results from \cref{section:uniformEO}.

Let $T$ be
an exponentially stable strongly continuous semigroup on a Hilbert
space $X$, $A$ its generator, $X_1$ the domain of $A$ with the norm
$\|v\|_1\coloneq \|Av\|$, and $C$ an
admissible observation operator taking values in a Hilbert space $Y$.
Assume that $A$ is normal and that $-A$ is sectorial of angle $<\pi
/2$. This is equivalent to $T$ being analytic and normal.
Let $\omega \in [0,\pi/2)$ such that $\sup_{z \in \sigma (A)}|
\arg(-z)|\leq \omega$.
Let $\Sigma \coloneq \{ z \in \CC_- \colon |\arg (-z)| \leq \omega \}$.
Throughout we make use of the spectral theorem for unbounded normal
operators as presented in~\cite{rudin_functional_1991}.
Let $E$ be the spectral measure of $A$. Assume that $(A,C)$
satisfies the Hautus test.
Let $m_h \in (0,\infty)$ such that for all $z \in \CC_-$ and $v \in X_1$ we have
\begin{equation*}
  m_h \| v \|^2 \leq \frac{1}{|\re z|^2}\| (z\id -A)v\|^2 +
  \frac{1}{|\re z|}\|C v\|^2.
\end{equation*}
Let $M \in (0,\infty)$ such that for all $v \in X_1$ we have
\begin{equation*}
  \int_0^\infty \| C T(s) v \|^2 \diff \lambda (s) \leq M \|v\|^2.
\end{equation*}
Such $M$ exists, because $C$ is admissible.

The Hautus test and admissibility have a scaling symmetry when
$(A,C)$ are replaced by $(t A, \sqrt{t}C)$ for $t \in (0,\infty)$ as the
following lemma shows.
\begin{lemma}\label{sector:scalingLemma}
  Let $t\in (0,\infty)$. Then for all $v \in X_1$ and $z \in \CC_-$ we have
  \begin{equation*}
    m_h \| v \|^2 \leq \frac{1}{|\re z|^2}\| (z \id - t A)v\|^2 + \frac{1}{|\re
    z|}\|\sqrt{t} C v\|^2
  \end{equation*}
  and for all $v \in X_1$ it is true that
  \begin{equation*}
    \int_0^\infty \| \sqrt{t} C  T( t s) v \|^2 \diff \lambda
    (s) = \int_0^\infty \|
    C T(s) v \|^2 \diff \lambda (s).
  \end{equation*}
\end{lemma}
\begin{proof}
  The first assertion follows from the fact that $t \cdot
  \mathbb{C}_- = \mathbb{C}_-$.
  The second assertion follows from the substitution $s \mapsto t s$.
\end{proof}
Now we choose the annulus used in the proof of
\cref{introduction:mainThm}. The specific choice for the radii will
be used in the proof of \cref{sector:dncnlim}.
Let $a^\prime , b^\prime \in (0,\infty)$ with $a^\prime <b^\prime$.
Let $q \in (1,\infty)$ such that
\begin{equation}
  \label{sector:qProperties}
  \frac{2}{q} <  a^\prime, \quad b^\prime < \log(q) , \quad
  \frac{2}{q} < e^{-b^\prime} (1-
  e^{-a^\prime}).
\end{equation}
Let  $a \coloneq\frac{2}{q}$ and $b\coloneq \frac{\log(q)}{\cos(\omega)}$.
For $t \in (0,\infty)$ let
\begin{equation*}
  \begin{split}
    \Omega_t \coloneq \{ z \in \Sigma \colon a\leq |z|t \leq b \}, \ P_t\coloneq
    E(\Omega_t), \ V_t \coloneq \ran P_t, \ C_t \coloneq \sqrt{t}
    C|_{V_t}, \ A_t
    \coloneq t A|_{V_t}.
  \end{split}
\end{equation*}
Let $t \in (0,\infty)$. We have that $V_t \subset
X_1$ with continuous embedding, because $\Omega_t$ is bounded.
This shows that $C_t$ is continuous.
Furthermore $\sigma (A_t) \subset t (\sigma (A) \cap \Omega_t) = \sigma (tA)
\cap \Omega_1\subset \Omega_1$.
Note also that $A_t$ generates the semigroup $T(\smallbullet t)|_{V_t}$.
In total this shows that ${(A_t)}_{t\in (0,\infty)}$ is uniformly
bounded and that ${(\exp(A_t \smallbullet))}_{t\in (0,\infty)}$ is
uniformly exponentially stable.
\cref{sector:scalingLemma} implies that ${(C_t)}_{t\in (0,\infty)}$
is uniformly admissible and that
${(A_t)}_{t\in (0,\infty)}$ and ${(C_t)}_{t\in (0,\infty)}$ satisfy
the Hautus test uniformly.
Therefore \cref{uniformEO:mainThm} gives constants $m_o \in
(0,\infty)$ and $L \in (0,\infty)$ such that
for all $t \in
(0,\infty)$ and $v \in V_t$
it is true that
\begin{equation}
  \label{sector:uniformEO}
  m_o \|v\|^2 \leq t \int_0^L \| C T(t s) v\|^2 \diff \lambda (s).
\end{equation}

Having obtained these local estimates we now consider how to glue
them together using powers of the following function.
Let $J \coloneq
(0,\infty) \setminus (a,b)$ and $f: \CC_- \to \CC$,
\begin{equation*}
  f(z) \coloneq e^z (1-e^z).
\end{equation*}
For $n \in \NN_+$ let $f_n : \CC_- \to \CC$,
$f_n (z)  \coloneq {f(z)}^n$ and
\begin{equation*}
  c_n \coloneq \inf_{\theta \in [-\omega,\omega]} \int_a^b \frac{1}{s}
  |f_n(-s e^{i\theta})|^2 \diff \lambda (s),
\end{equation*}
\begin{equation*}
  d_n \coloneq \sup_{\theta \in [-\omega,\omega]} \int_{J} \frac{1}{s}
  |f_n(-s e^{i\theta})|^2 \diff \lambda (s).
\end{equation*}
\begin{proposition}
  \label{sector:dncnlim}
  For all $n \in \NN_+$ we have $|f_n|$ is bounded on $\Sigma$, $0< c_n
  <\infty$ and $d_n <\infty$.
  Moreover, it is true that
  \begin{equation*}
    \lim_{n \to \infty} \frac{d_n}{c_n} =0.
  \end{equation*}
\end{proposition}
\begin{proof}
  Note that $\cos (\omega ) >0$, because $\omega \in
  [0,\pi/2)$.
  Let $s \in (0,\infty)$ and $\theta \in [-\omega, \omega]$.
  Then by the reverse triangle inequality we have
  \begin{equation}
    \label{sector:fEstimBelow}
    |f(-s e^{i\theta}) | \geq e^{-s \cos (\theta)} (1- e^{-s \cos
    (\theta)})
  \end{equation}
  and by the triangle inequality
  \begin{equation}
    \label{sector:fEstimAboveInfty}
    | f(-s e^{i\theta})| \leq 2 e^{-s \cos (\omega)}.
  \end{equation}
  Furthermore we get
  \begin{equation}
    \label{sector:fEstimAboveZero}
    | f(-s e^{i\theta})| \leq s,
  \end{equation}
  because the fundamental theorem of calculus shows that
  \begin{equation*}
    1-\exp(- s e^{i\theta})  =  s e^{i \theta} \int_0^1 \exp (- t s
    e^{i \theta})
    \diff \lambda (t)
  \end{equation*}
  from which we obtain that
  \begin{equation*}
    |1-\exp(- s e^{i\theta})| \leq s \int_0^1 \underbrace{\exp ( - ts
    \cos (\theta) )}_{\leq 1} \diff \lambda (t) \leq s.
  \end{equation*}
  In total, this shows that for all $n \in \NN_+$ we have that
  $|f_n|$ is bounded on
  $\Sigma$, $0< c_n
  <\infty$, and $d_n <\infty$.

  Let $n \in \NN_+$ and $\theta \in [-\omega,\omega]$.
  By~\eqref{sector:qProperties} we have $b\geq b^\prime/\cos(\theta)$
  and $a\leq a^\prime$. This together with~\eqref{sector:fEstimBelow}
  implies that
  \begin{equation*}
    \begin{split}
      \int_a^b s^{-1} | f_n(-s e^{i\theta})|^2 \diff \lambda (s) &\geq
      \int_{a^\prime /\cos(\theta)}^{b^\prime /\cos(\theta)} s^{-1} |
      f_n(-s e^{i\theta})|^2 \diff \lambda (s) \\
      &\geq
      {\big(\underbrace{e^{-b^\prime} (1- e^{-a^\prime})}_{\eqcolon
      p} \big)}^{2n} \log( b^\prime/a^\prime).
    \end{split}
  \end{equation*}
  Therefore we obtain
  \begin{equation*}
    c_n \geq  {p}^{2n} \log( b^\prime/a^\prime) >0.
  \end{equation*}
  On the other hand from~\eqref{sector:fEstimAboveZero} we get
  \begin{equation*}
    \int_0^a s^{-1} |f_n (-s e^{i\theta}) |^{2} \diff \lambda (s)
    \leq \frac{a^{2n}}{2n}
  \end{equation*}
  and
  \begin{equation*}
    \int_b^\infty s^{-1} |f_n (-s e^{i\theta}) |^{2} \diff \lambda (s)
    \leq  b^{-1} 4^n \int_b^\infty e^{- 2 n s \cos (\omega)} \diff
    \lambda (s) =
    \frac{{(2 e^{- b \cos (\omega)})}^{2n}}{  2n \cos(\omega) b } .
  \end{equation*}
  Therefore, using $2e^{- b\cos (\omega)} = 2/q =  a$, we have
  \begin{equation*}
    d_n \leq \frac{ 1+ {(\cos(\omega)b )}^{-1}}{2 n} {a}^{2n}.
  \end{equation*}
  This shows that for all $n \in \NN_+$ it is true that
  \begin{equation*}
    \frac{d_n}{c_n} \leq \frac{1}{\log(b^\prime/a^\prime)} \frac{ 1+
    {(\cos(\omega)b )}^{-1}}{2 n} {\bigg( \frac{a}{p}\bigg)}^{2n}.
  \end{equation*}
  Therefore $\lim_{n\to \infty} d_n/ c_n =0$, since $a/p<1$
  by~\eqref{sector:qProperties}.
\end{proof}
In the following we will use the next lemma, without mention, to
conclude the measurability of certain integrands. For example we will
apply it to $g(t,z) \coloneq f_n (tz)\mathbf{1}_{\Omega_t} (z)$.
\begin{lemma}\label{sector:measurabilityLemma}
  Let $v\in X_1$ and $ g : (0,\infty) \times \CC_- \to \CC$
  measurable, such that for every $t \in (0,\infty)$ the map
  $g(t, \smallbullet)$ is bounded on $\sigma(A)$.
  Then $(0,\infty)\ni t \mapsto g (t, A) v \in X_1$ is strongly measurable.
  Furthermore $(0,\infty) \ni t  \mapsto \int_0^L \| C T(st) g(t,A)v \|^2
  \diff \lambda (s)$ is measurable.
\end{lemma}
\begin{proof}
  Let $t \in (0,\infty)$. Note that $g(t,A) v \in X_1$,
  because $g(t,\smallbullet)$ is bounded on $\sigma (A)$ and the
  multiplication theorem for the functional calculus.
  The map $U : L^2 (\CC_-,|\smallbullet|^2\|
  E(\smallbullet)  v \|^2) \to X_1$, $U f \coloneq \int_{\CC_-} f (z)
  \diff E (z) v$ is an
  isometry and for all $t \in
  (0,\infty)$ we have $g (t,A) v \in \ran U$.
  The subspace $\ran U$ is closed, since $U$ is an isometry, and
  separable, because the involved measure
  is finite (since $v \in X_1$) and the Borel $\sigma$-algebra of $\CC_-$ is
  countably generated together with a standard result from measure
  theory~\cite[Proposition 3.4.5,
  p.~102]{cohn_measure_2013}. Therefore we can apply the Pettis
  measurability theorem.
  It is left to show that $t \mapsto g(t,A)v$ is weakly measurable.
  Let $ w \in \ran U$. Let $h \in L^2 (\CC_-,|\smallbullet|^2\|
  E(\smallbullet)  v \|^2)$ such that $w = U h$.
  Then
  \begin{equation*}
    \langle g(t,A) v , w\rangle_{X_1} = \int_{\CC_-}  g(t,z) \bar h(z)
    |z|^2\diff \|E(z) v \|^2.
  \end{equation*}
  The integrand is measurable on $(0,\infty) \times \CC_-$ and the
  involved measure is $\sigma$-finite. Therefore
  part of Tonelli's theorem shows that $t \mapsto \langle g(t,A) v ,
  w\rangle$ is measurable.
  Since $w \in \ran U$ is arbitrary we get that $t \mapsto g(t,A) v$
  is weakly measurable and hence, by the
  Pettis measurability theorem~\cite[Theorem 1.1.6,
  p.~5]{hytonen_analysis_2016}, that it is strongly measurable.

  For the second assertion: $T|_{X_1}$ is strongly continuous.
  Therefore $(0,\infty) \times (0,L) \ni (t,s) \mapsto C T(st)|_{X_1}$
  is also strongly continuous.
  By what has been shown $(0,\infty) \times (0,L) \ni (t,s) \mapsto
  g(t,A)v$ is strongly measurable.
  Therefore~\cite[Proposition 1.1.28, p.~12]{hytonen_analysis_2016}
  implies that $(0,\infty) \times (0,L) \ni  (t,s) \mapsto
  C T(st) g(t,A)v$ is strongly measurable.
  The conclusion follows, again, from Tonelli's theorem.
\end{proof}
Using the spectral theorem we obtain the following two estimates.
The first yields~\eqref{eq:intro1}, and the second will be used later
together with the admissibility of $C$ to obtain~\eqref{eq:intro3}.
\begin{proposition}
  \label{sector:cnDnFormula}
  Let $n \in \NN_+$ and $v \in X_1$. Then we have
  \begin{equation*}
    \int_0^\infty t^{-1}\| f_n (t A) P_t v\|^2 \diff \lambda (t)
    \geq c_n \|v\|^2
  \end{equation*}
  and
  \begin{equation*}
    \int_0^\infty t^{-1}\| f_n(t A) (I-P_t) v\|^2 \diff \lambda (t)
    \leq d_n \|v\|^2.
  \end{equation*}
\end{proposition}
\begin{proof}
  Let $ v \in X_1$ and $ n \in \NN_+$.
  Then using Tonelli's theorem and the spectral theorem gives
  \begin{equation*}
    \begin{split}
      \int_0^\infty t^{-1}\| f_n (t A) P_t v\|^2 \diff \lambda (t) =
      \int_{\Sigma} \int_0^\infty  t^{-1} | f_n(t z)|^2
      \mathbf{1}_{[a,b]}(t |z|) \diff \lambda (t) \diff \| E(z) v\|^2.
    \end{split}
  \end{equation*}
  Let $z \in \Sigma$. Let $\theta \in [-\omega, \omega]$ (the unique
  one) such that $z= -|z| e^{i \theta}$.
  Then using the substitution $t \mapsto |z| t$ shows that
  \begin{equation*}
    \int_0^\infty  t^{-1} | f_n(t z)|^2
    \mathbf{1}_{[a,b]}(t |z|) \diff \lambda (t) =
    \int_a^b  s^{-1} | f_n(-s e^{i\theta})|^2  \diff \lambda (s) \geq c_n.
  \end{equation*}
  Therefore
  \begin{equation*}
    \int_0^\infty t^{-1}\| f_n (t A) P_t v\|^2 \diff \lambda (t) \geq
    c_n \int_{\Sigma}  \diff \| E(z) v\|^2= c_n \|v\|^2.
  \end{equation*}
  The proof of the second inequality runs along the same lines.
\end{proof}
To obtain the estimate~\eqref{eq:intro2} we use the fact that $f_n$
is a linear combination of functions of the form $z \mapsto \exp
(\beta z)$, $\beta >0$ for all $n \in \NN_+$.
\begin{lemma}\label{sector:dblIntEstim}
  Let $N \in \NN_+$, $\alpha \in \CC^N$, $\alpha \neq 0$, $\beta \in
  {(0,\infty)}^N$
  and $h : \CC_- \to \CC$, $h(z) \coloneq \sum_{j=1}^N \alpha_j \exp
  (\beta_j z)$.
  Then there exists $\kappa \in (0,\infty)$ such that for all $v \in
  X_1$ we have
  \begin{equation*}
    \int_0^\infty \int_0^L  \| C T(s t) h( t A) v \|^2 \diff \lambda
    (s) \diff \lambda (t)  \leq \kappa \int_0^\infty \| C T(t) v\|^2
    \diff \lambda (t).
  \end{equation*}
\end{lemma}
\begin{proof}
  Let $v \in X_1$ and
  \begin{equation*}
    I_1 \coloneq \int_0^\infty \int_0^L  \| C T(s t) h( t A) v \|^2
    \diff \lambda
    (s) \diff \lambda (t) .
  \end{equation*}
  From the triangle inequality, Jensen's inequality (applied to the
  convex function $t \mapsto t^2$), and the multiplication theorem for
  the functional calculus we get
  \begin{equation*}
    \begin{split}
      I_1 &\leq N \sum_{j=1}^N {|\alpha_j|}^2
      \underbrace{ \int_0^\infty \int_0^L  \| C \exp (t (s+ \beta_j)
        A) v \|^2 \diff \lambda
      (s) \diff \lambda (t)}_{\eqcolon I_{2,j}} .
    \end{split}
  \end{equation*}
  Let $j \in \{1,\dots, N\}$.
  Tonelli's theorem and the substitution $t \mapsto t (s+
  \beta_j)$ imply that
  \begin{equation*}
    I_{2,j} = \int_0^L  \frac{1}{s+ \beta_j} \diff \lambda (s)
    \int_0^\infty \| C T(t)
    v\|^2 \diff \lambda (t).
  \end{equation*}
  The fundamental theorem of calculus yields
  \begin{equation*}
    \int_0^L  \frac{1}{s+ \beta_j} \diff \lambda (s) = \log ( L+
    \beta_j) - \log ( \beta_j) = \log \bigg(1+\frac{L}{\beta_j}\bigg)>0.
  \end{equation*}
  And so in total we have
  \begin{equation*}
    I_1 \leq \underbrace{N \sum_{j=1}^N {|\alpha_j|}^2 \log
    \bigg(1+\frac{L}{\beta_j}\bigg)}_{>0} \int_0^\infty \| C T(t)
    v\|^2 \diff \lambda (t). \qedhere
  \end{equation*}
\end{proof}
\begin{proof}[\textbf{\upshape{Proof of \cref{introduction:mainThm}}}]
  Let $v \in X_1$ and $n \in \NN_+$. Then from~\eqref{sector:uniformEO}
  and \cref{sector:cnDnFormula} we get
  \begin{equation}
    \label{sector:EOEstimate1}
    \begin{split}
      m_o c_n \|v\|^2 &\leq m_o \int_0^\infty t^{-1}\| f_n(t A) P_t
      v\|^2 \diff
      \lambda (t)  \\
      &\leq   \int_0^\infty \int_0^L  \|C T(ts) f_n(t A) P_t v\|^2
      \diff \lambda (s) \diff
      \lambda (t)  \\
      &\leq 2 \bigg(\underbrace{\int_0^\infty \int_0^L  \|C T(t s)
          f_n(t A) v\|^2
          \diff \lambda (s) \diff
        \lambda (t)}_{\eqcolon I_{1,n,v}} \\ &\mathrel{\phantom{=}}+
        \underbrace{\int_0^\infty \int_0^L
          \| C T(t s) f_n(t A)
          (I- P_t) v\|^2
          \diff \lambda (s) \diff
      \lambda (t)}_{\eqcolon I_{2,n,v}} \bigg).
    \end{split}
  \end{equation}
  Now, using the substitution $s\mapsto  t s$ for $t \in
  (0,\infty)$ in the inner integral and extending the upper integral
  bound to $\infty$ it follows that for all $n \in \NN_+$ and $v \in
  X_1$ we have
  \begin{equation*}
    \begin{split}
      I_{2,n, v}&\leq
      \int_0^\infty \int_0^{\infty} \frac{1}{t} \| C T(s) f_n(t A)
      (I- P_t) v\|^2
      \diff \lambda (s) \diff
      \lambda (t) \\
      &\leq M d_n \|v\|^2.
    \end{split}
  \end{equation*}
  Let $N \in \NN_+$ such that $m_o c_N - 2 M d_N >0$. The
  number $N$ exists
  by \cref{sector:dncnlim}.
  Let $\kappa \in (0,\infty)$ such that for all $v \in X_1$:
  \begin{equation*}
    I_{1,N,v} = \int_0^\infty \int_0^L  \|C T(t s) f_N(t A) v\|^2
    \diff \lambda (s) \diff
    \lambda (t) \leq \kappa \int_0^\infty  \| C T(t) v\|^2 \diff
    \lambda (t).
  \end{equation*}
  Such $\kappa$ exists by \cref{sector:dblIntEstim}.
  Then~\eqref{sector:EOEstimate1} implies that for all $v \in X_1$ it
  is true that
  \begin{equation*}
    \underbrace{(m_o c_N - 2 M d_N)}_{>0} \|v\|^2 \leq 2 \kappa
    \int_0^\infty  \| C T(t)
    v\|^2 \diff \lambda (t). \qedhere
  \end{equation*}
\end{proof}

\section{Criteria for applications}\label{sec:criteria}
The examples in \cref{sec:example1} and \cref{sec:example2} illustrate how
\cref{introduction:mainThm} can establish exact observability even
when the lower
comparison with the instantaneous energy dissipation required by
\cref{examples:RealPartEOSufficiencyEstimate} fails.

In the first subsection we present the energy dissipation comparison
criterion for checking
admissibility and exact observability.

In the second subsection we give a sufficient condition for verifying
the Hautus test
using only the real part of the generator.
\subsection*{Energy dissipation criteria for admissibility and exact
observability}
Let $A$ be a normal generator of an exponentially stable strongly
continuous semigroup $T$ on the Hilbert space $X$. In this setting
there is always a
natural choice for the
observation operator as the following energy identity shows.
Define the \emph{energy} $\mathscr{E} : X \to [0,\infty)$ by
$\mathscr{E}(v) \coloneq
\frac{1}{2}\|v\|^2$.
Then for every $v \in X_1$ and $t \in [0,\infty)$ we have
\begin{equation*}
  \begin{split}
    (\mathscr{E} {(T(\smallbullet) v))}^\prime (t) &= \re \langle A
    T(t) v, T(t)v\rangle \\
    &= \langle \re A T(t) v, T(t)v\rangle \\
    &= - \| {(- \re A)}^{1/2} T(t)v\|^2.
  \end{split}
\end{equation*}
Therefore $\| {(- \re A)}^{1/2} \smallbullet \|^2$ gives the
instantaneous rate of
energy dissipation along the
orbits of $T$.
By integration and exponential stability of $T$ we obtain for all $v
\in X_1$ that
\begin{equation}
  \label{criteria:EnergyIdentity}
  \int_0^\infty \| {(- \re A)}^{1/2} T(s) v\|^2 \diff \lambda(s) =
  \mathscr{E}(v)=
  \frac{1}{2} \| v\|^2.
\end{equation}
The operator ${(- \re A)}^{1/2}:  X_1 \to X$ is continuous since for
all $v \in X_1$ we have
\begin{equation*}
  \| {(- \re A)}^{1/2} v\|^2 = \re \langle - A v, v \rangle \leq \|A
  v\| \| v\| = \| Av
  \| \|A^{-1} A v \|\leq \|A v \|^2  \|A^{-1}\|.
\end{equation*}
This shows together with~\eqref{criteria:EnergyIdentity} that ${(-
\re A)}^{1/2}$ is
admissible and exactly observable.

From these considerations we obtain the following energy dissipation
criteria for
admissibility and exact
observability for an observation operator $C$ by comparing $\|C
\smallbullet\|^2$ to the
instantaneous energy
decay rate given by $\| {(- \re A)}^{1/2} \smallbullet \|^2$.
\begin{lemma}
  \label{examples:RealPartAdmSufficiencyEstimate}
  Let $C : X_1 \to Y$ be bounded. If there exists $c  \in (0,\infty)$
  such that for all $v \in X_1$ it is true that
  \begin{equation*}
    c \| C v\| \leq  \| {(- \re A)}^{1/2} v\|
  \end{equation*}
  then $C$ is admissible.
\end{lemma}
\begin{proof}
  This follows from the fact that ${(- \re A)}^{1/2}$ is admissible.
\end{proof}
\begin{lemma}
  \label{examples:RealPartEOSufficiencyEstimate}
  Let $C : X_1 \to Y$ be bounded and admissible. If there exists $c
  \in (0,\infty)$ such
  that for all $v \in X_1$ it is true that
  \begin{equation*}
    \| {(- \re A)}^{1/2} v\| \leq c \| C v\|
  \end{equation*}
  then $(A,C)$ is exactly observable.
\end{lemma}
\begin{proof}
  This follows from the fact that ${(- \re A)}^{1/2}$ is exactly observable.
\end{proof}

Using \cref{examples:RealPartEOSufficiencyEstimate} is a simple and
often successful way of
proving exact observability.
In our applications of \cref{introduction:mainThm} in \cref{sec:example1} and
\cref{sec:example2} we therefore consider
examples where the generated semigroup $T$ is analytic and the assumptions of
\cref{examples:RealPartEOSufficiencyEstimate} are not
satisfied, but the Hautus test is. In both of our examples $A$ is
unbounded, which, if
$T$ is analytic, is equivalent to $\inf_{z \in \sigma (A)} \re z = -\infty$.
\subsection*{Real-part criterion for the Hautus test}
The following lemma shows that verifying a Hautus test type inequality for the
real part and real parameters is sufficient to conclude that the
Hautus test is valid for normal operators.
It will be used to verify the Hautus test for our examples in
\cref{sec:example1} and
\cref{sec:example2}.
\begin{lemma}\label{examples:lemRealHTImpliesHT}
  Let $A$ be a normal generator of an exponentially stable strongly
  continuous semigroup and $C$ an admissible observation operator.
  Assume that there exists $m \in (0,\infty)$ such that for all
  $s \in (0,\infty)$ and $v \in X_1$ it is true that
  \begin{equation*}
    m \| v \|^2 \leq s^{-2} \| (s+ \re A) v\|^2 + s^{-1} \| C v\|^2
  \end{equation*}
  Then $(A,C)$ satisfies
  the Hautus test.
\end{lemma}
\begin{proof}
  Let $z \in \CC_-$ and  $v \in X_1$. Then we have
  \begin{equation*}
    \begin{split}
      m \| v \|^2 &\leq \frac{1}{{| \re z|}^2} \| (|\re z| \id + \re A) v\|^2 +
      \frac{1}{|\re z|} \| C v\|^2 \\
      &\leq \frac{1}{{| \re z|}^2} \big (\|  (\re z \id - \re A)
      v\|^2 +  \| (\im z \id - \im A) v\|^2 \big)
      +\frac{1}{|\re z|} \| C v\|^2 \\
      &= \frac{1}{{| \re z|}^2} \| (z \id -  A)
      v\|^2 +
      \frac{1}{|\re z|} \| C v\|^2,
    \end{split}
  \end{equation*}
  where for the last equality it is used that $A$ is normal.
\end{proof}
\begin{remark}
  The reverse implication in \cref{examples:lemRealHTImpliesHT} is false.
  A counterexample is given by $X \coloneq \CC^2$, $Y \coloneq \CC$ and
  \begin{equation*}
    A \coloneq
    \begin{pmatrix}
      -1+i & 0 \\
      0 & -1-i
    \end{pmatrix}, \quad
    C \coloneq
    \begin{pmatrix}
      1 & 1
    \end{pmatrix}.
  \end{equation*}
  Then $(A,C)$ is exactly observable, as any of the usual
  finite-dimensional tests
  quickly show, and therefore $(A,C)$ satisfies the
  Hautus test.
  However $\re A = -\id$, which implies $ \id + \re A = 0$ and $C
  {(-1, 1)}^\top = 0$.
\end{remark}

\section{Exact observability of a multiplication semigroup from high-pass frequency domain
measurements of the weighted state}\label{sec:example1}
In this section we give an example application of \cref{introduction:mainThm} to a
multiplication operator with continuous spectrum.
In this example the assumptions of the energy dissipation criterion for exact observability,
\cref{examples:RealPartEOSufficiencyEstimate}, are not
satisfied (see \cref{example1:EOlemmaNotSatisfied}).

Let $X\coloneq L^2((0,\infty))$, $Y\coloneq L^2(\RR)$, and
$\mathscr{F}$ be the unitary Fourier transform.
Throughout this section we will sometimes implicitly view elements of
$X$ as elements of
$L^2(\RR)$ by zero extension.
Let $\nu \in L^\infty ((0,\infty))$ be real-valued almost everywhere.
Define a normal operator $A$ on $X$ by
\begin{equation*}
  (A f)(\xi) \coloneq - e^\xi (1+ i \nu (\xi)) f(\xi)
\end{equation*}
with the maximal domain.
Let $\delta \in (0,\infty)$ and define $C : X_1 \to Y$ by
\begin{equation*}
  C f \coloneq \mathbf{1}_{\mathbb{R} \setminus [-\delta, \delta]}
  \mathscr{F}  {(-\re A)}^{1/2}  f.
\end{equation*}
It is clear that $A$ is normal and that $A$ generates an
exponentially stable analytic strongly continuous semigroup.
The boundedness of $C :X_1 \to Y$ is also evident.
For $f \in X_1$ the observation $C f$ is the truncation of the frequency band $[-\delta,
\delta]$ from the Fourier transform of the natural
observation ${(-\re A)}^{1/2} f$, which in this case is just $f$ weighted with
$e^{\smallbullet /2}$.
We will show that the removal of the finite frequency band $[-\delta, \delta]$ does
not destroy the exact
observability using the Hautus test.
\begin{lemma}
  The observation operator $C$ is admissible.
\end{lemma}
\begin{proof}
  For every $f \in X_1$ we have $\|Cf\| \leq \| {(- \re A)}^{1/2} f\|$.
  Therefore $C$ is admissible by
  \cref{examples:RealPartAdmSufficiencyEstimate}.
\end{proof}
To show that $(A,C)$ satisfies the Hautus test we need the following lemma.
\begin{lemma}\label{examples:FourierOfCompactSupp}
  Let $r_0, r_1 \in \RR$ with $r_0< r_1$ and let $h \in X$ vanish outside
  of $[r_0,r_1]$ almost
  everywhere. Let $\omega \in \RR$. Then we have
  \begin{equation*}
    | (\mathscr{F} h)(\omega)| \leq \|h \| \sqrt{\frac{r_1-r_0}{2 \pi }}.
  \end{equation*}
\end{lemma}
\begin{proof}
  Using the Cauchy-Schwarz inequality we get
  \begin{equation*}
    \begin{split}
      | (\mathscr{F} h)(\omega)|
      &=
      \frac{1}{\sqrt{2\pi}} \bigg|\int_{r_0}^{r_1}
      e^{i \xi \omega } h(\xi) \diff \lambda (\xi)\bigg| \\
      &\leq \frac{1}{\sqrt{2\pi}} \int_{r_0}^{r_1}|h(\xi)| \diff
      \lambda (\xi)  \\
      &\leq\frac{1}{\sqrt{2\pi}} \sqrt{r_1-r_0} \|h\|. \qedhere
    \end{split}
  \end{equation*}
\end{proof}
\begin{proposition}
  The system $(A,C)$ satisfies the Hautus test and so, by
  \cref{introduction:mainThm}, $(A,C)$ is exactly observable.
\end{proposition}
\begin{proof}
  To verify the Hautus test we will use \cref{examples:lemRealHTImpliesHT}.
  To this end let $s \in (0,\infty)$.
  Let $\alpha, \beta \in (0,\infty)$ such that $\alpha <1 < \beta$ and
  \begin{equation*}
    \frac{\delta}{\pi} \log \bigg(\frac{\beta}{\alpha}\bigg) <1.
  \end{equation*}
  Let $r_1 \coloneq \log( s \beta)$ and $r_0 \coloneq \log(s \alpha)$.
  Note that $r_1 -r_0 = \log( \beta / \alpha) >0$.
  Let $J \coloneq (0,\infty) \cap [r_0, r_1]$.
  Let $f \in X_1$, $f_1 \coloneq \mathbf{1}_J f$, $f_2 \coloneq
  \mathbf{1}_{\RR \setminus J} f$.
  Then by the spectral theorem
  \begin{equation*}
    \| {(- \re A)}^{1/2}f_1\|^2 \geq e^{r_0} \|f_1\|^2 = s
    \alpha \|f_1\|^2.
  \end{equation*}
  Furthermore
  \begin{equation*}
    \begin{split}
      \| {(- \re A)}^{1/2}f_1\|^2  &=\|  \mathscr{F} {(- \re
      A)}^{1/2}f_1\|^2  \\
      &\leq \|  C f_1\|^2 +  \|
      \mathbf{1}_{[-\delta,\delta]} \mathscr{F} {(- \re
      A)}^{1/2}f_1\|^2  \\
      &\leq \|  C f_1\|^2 + \delta \frac{r_1-r_0}{\pi}
      \| {(- \re A)}^{1/2}f_1\|^2,
    \end{split}
  \end{equation*}
  where \cref{examples:FourierOfCompactSupp} is used for the last inequality.
  In total we have shown that
  \begin{equation*}
    s \underbrace{\alpha \bigg(1-  \frac{\delta}{\pi} \log \bigg
    (\frac{\beta}{\alpha}\bigg) \bigg)}_{\eqcolon \kappa} \|f_1\|^2
    \leq \| C f_1 \|^2.
  \end{equation*}
  Moreover, it is true that
  \begin{equation*}
    \begin{split}
      \| C f_2\|^2 &\leq \| {(-\re A)}^{1/2} f_2\|^2 \\
      &\leq   \sup_{\xi \in
      \mathbb{R} \setminus [r_0,r_1]}\frac{e^\xi}{{(s  - e^\xi)}^2}
      \| (s+ \re A) f_2 \|^2\\
      &\leq \frac{\beta}{s \underbrace{\min\{ (\beta-1)^2,
      (1-\alpha)^2 \}}_{\eqcolon \gamma }} \| (s+ \re A) f \|^2.
    \end{split}
  \end{equation*}
  By the same argument we also get
  \begin{equation*}
    \|f_2\|^2 \leq \frac{1}{s^2 \gamma } \| (s + \re A) f \|^2.
  \end{equation*}
  Since $\| C f_1\|^2 \leq 2 (\|C f\|^2 +
  \|Cf_2\|^2)$, we obtain
  \begin{equation*}
    \begin{split}
      \| f \|^2 &= \| f_1\|^2 + \|f_2\|^2 \\
      &\leq \kappa^{-1}\frac{1}{s}\| Cf_1\|^2 + \gamma^{-1}
      \frac{1}{s^2} \| (s \id + \re A) f\|^2 \\
      &\leq 2 \kappa^{-1} \frac{1}{s}  \| Cf \|^2 + \gamma^{-1} ( 2
      \kappa^{-1} \beta + 1) \frac{1}{s^2} \| (s \id + \re A) f\|^2.
    \end{split}
  \end{equation*}
  In view of \cref{examples:lemRealHTImpliesHT} this shows that
  $(A,C)$ satisfies the Hautus test.
\end{proof}
Finally we show that the assumptions of \cref{examples:RealPartEOSufficiencyEstimate} are
not satisfied.
\begin{proposition}
  \label{example1:EOlemmaNotSatisfied}
  There does not exist $c \in (0,\infty)$ such
  that for all $v \in X_1$ it is true that
  \begin{equation*}
    \| {(- \re A)}^{1/2} v\| \leq c \| C v\|.
  \end{equation*}
\end{proposition}
\begin{proof}
  Let $g \in C_c((0,\infty))$ with $\|g\|_X=1$. Define a sequence
  ${(g_n)}_{ n \in \NN_+}$ in $X$ by
  \begin{equation*}
    g_n (\xi) \coloneq \frac{1}{\sqrt{n}} g( n^{-1} \xi).
  \end{equation*}
  Then $\|g_n\| =1$ for all $n \in \NN_+$.
  Now define a sequence
  ${(f_n)}_{ n \in \NN_+}$ in $X$ by
  \begin{equation*}
    f_n \coloneq {(- \re A)}^{-1/2} g_n.
  \end{equation*}
  Let $n \in \NN_+$. Clearly $\|{(-\re A)}^{1/2}f_n\|= \| g_n \| = 1$.
  On the other hand by definition and the substitution rule
  \begin{equation*}
    {\| C f_n\|}^2 = \int_{\RR \setminus n [-\delta, \delta]} |
    (\mathscr{F} g)(\omega) |^2 \diff  \lambda (\omega).
  \end{equation*}
  Since $\mathscr{F}g \in L^2(\RR)$ we get from Lebesgue's theorem of
  dominated convergence that
  \begin{equation*}
    \lim_{n \to \infty} \|C f_n \|^2 =0. \qedhere
  \end{equation*}
\end{proof}

\section{Exact observability of derivative observation for a
  diffusion equation with a
singular potential}\label{sec:example2}
In this section we use \cref{introduction:mainThm} to show exact
observability for a one-dimensional
diffusion equation with a singular
potential and derivative observation.
This example also illustrates a strategy for verifying the Hautus test that
may apply to other diffusion equations (see
\cref{example2:generalizationRemark}).

In particular with our potential the assumptions of the energy
dissipation criterion for
exact observability,
\cref{examples:RealPartEOSufficiencyEstimate}, are not
satisfied (see \cref{examples2:No-ReAsqEstimate}).
Our proof does not rely on an explicit representation of the generated
semigroup or the spectral measure of the generator.
This highlights the viability of using the Hautus test together with
\cref{introduction:mainThm} to show exact observability.

Let $k, \ell \in \mathbb{N}$, $k \geq 3$, $k+2 \leq \ell \leq 2k$.
For $j \in \NN_+$ let $x_j \coloneq 2^{-j}$, $d_j \coloneq 2^{- k j}$,
$H_j \coloneq 2^{ \ell j}$.
Moreover, let $\psi \in C_c^\infty((-1,1))$ such that $\psi(0)=1$ and $0\leq
\psi \leq 1$.
\begin{lemma}
  Let $j \in \NN_+$. Then we have that $\supp \psi\big(
  \frac{\smallbullet -x_j}{d_j} \big) \subset [
  x_j-d_j ,x_j+d_j]$ and
  \begin{equation*}
    [x_j -d_j , x_j + d_j] \cap [x_{j+1} -d_{j+1},
    x_{j+1} + d_{j+1}]= \varnothing.
  \end{equation*}
\end{lemma}
\begin{proof}
  The fact that $\supp \psi\big(
  \frac{\smallbullet -x_j}{d_j} \big) \subset [
  x_j-d_j ,x_j+d_j]$ is clear.
  From the definition of $x_j$ and $d_j$ we get
  \begin{equation*}
    x_j - d_j = 2^{-j} - 2^{ -k j} = 2^{-j} (1- 2^{-(k-1) j})
  \end{equation*}
  and
  \begin{equation*}
    x_{j+1} + d_{j+1} = 2^{-(j+1)} + 2^{-k (j+1)} = 2^{-j} 2^{-1} (1 +
    2^{-(k-1) (j+1)}).
  \end{equation*}
  Therefore it is true that
  \begin{equation*}
    \begin{split}
      x_{j+1}+ d_{j+1} < x_j -d_j &\Leftrightarrow 2^{-1} (1 +
      2^{-(k-1) (j+1)}) < 1- 2^{-(k-1)j} \\
      &\Leftrightarrow
      2^{-(k-1)j } (1+ 2^{-k}) < 2^{-1}.
    \end{split}
  \end{equation*}
  Using $k \geq 3$ and $j \geq 1$ we obtain
  \begin{equation*}
    2^{-(k-1)j } (1+ 2^{-k})  \leq
    2^{-(k-1)} (1+ 2^{-k}) \leq 2^{-2} (1+ 2^{-3}) = \frac{1}{4}\cdot
    \frac{9}{8} = \frac{9}{32} <\frac{1}{2}. \qedhere
  \end{equation*}
\end{proof}
Define $V: (0,1) \to [0,\infty)$,
\begin{equation*}
  V(x) \coloneq \sum_{j=1}^\infty H_j \psi \bigg( \frac{x-x_j}{d_j} \bigg).
\end{equation*}
The series is well defined, because the sum is actually locally finite.
This also implies that $V \in C^\infty ((0,1))$.

Let $X \coloneq L^2((0,1))$ and $Y \coloneq L^2((0,1))$.
Let
\begin{equation*}
  \mathcal{D} \coloneq \{ f \in H_0^1 ((0,1))\colon \sqrt{V} f \in
  L^2 ((0,1)) \}
\end{equation*}
and define $\mathfrak{b} \colon \mathcal{D} \times \mathcal{D} \to \CC$,
\begin{equation*}
  \mathfrak{b}(f,g) \coloneq \int_0^1 f^\prime \bar g^\prime + f \bar g V
  \diff \lambda.
\end{equation*}
Then $\mathfrak{H} \coloneq (\mathcal{D}, \mathfrak{b})$ is a Hilbert
space (using
Poincar\' e's inequality) and
$\mathfrak{H}$ is densely and continuously embedded in $X$.
Let
\begin{equation*}
  \mathscr{D} \coloneq \{ f \in \mathcal{D} \colon \exists h \in X
  \ \forall g \in \mathcal{D} \ \langle h,g\rangle = \mathfrak{b}(f,g)\}.
\end{equation*}
For $f \in \mathcal{D}$ there exists at most one $h \in X$ such
that for all $g \in \mathcal{D}$ we have $\langle h,g\rangle =
\mathfrak{b}(f,g)$, because
$C_c^\infty ((0,1)) \subset \mathcal{D}$ is dense in $X$.
Therefore we can define a linear operator $B : \mathscr{D} \to X$,
$Bf \coloneq h$, where
$h \in X$ is the unique one such that for all $g \in \mathcal{D}$ we
have $\langle h,g\rangle = \mathfrak{b}(f,g)$.
By~\cite[Chapter VI, Theorems 2.1 and 2.6,
pp.~322--323]{kato_perturbation_1995} the
operator $B$ is self-adjoint.
It is easily seen, by using integration by parts and the fact that $V
\in C^\infty
((0,1))$, that $C_c^\infty ((0,1)) \subset \dom
(B)$ and that for all $f \in
C_c^\infty ((0,1))$ we have
$B f= - f^{\prime \prime} + V f$.

Let $A \coloneq - B$. For all $f \in \dom (A)$ we have
\begin{equation}
  \label{example:SLInequality1}
  \begin{split}
    \langle - A f, f\rangle &= \mathfrak{b}(f,f)  \geq
    \|f^\prime \|^2 \geq \|f\|^2,
  \end{split}
\end{equation}
where we use
Poincar\'e's inequality and the fact that $V \geq 0$.
This shows that $\sigma (A) \subset (-\infty,-1]$.
Since $A$ is unbounded we also have $\inf_{z\in \sigma(A)} \re z = -\infty$.
Therefore $A$ generates an exponentially stable strongly continuous semigroup.
Let $C : X_1 \to Y$, $C f= f^\prime$.
Let $f \in \dom (A)$. Then by~\eqref{example:SLInequality1} it is true that
\begin{equation*}
  \| C f \|^2 \leq \langle - Af,f \rangle \leq \| Af \| \| f\| =  \|A
  f\| \|   A^{-1} A f\| \leq \| A^{-1}\| \| A f\|^2.
\end{equation*}
Therefore $C$ is continuous.
Furthermore~\eqref{example:SLInequality1} implies
\begin{equation*}
  \| C f\| \leq  \| {(- A)}^{1/2} f  \| .
\end{equation*}
Therefore \cref{examples:RealPartAdmSufficiencyEstimate} implies that
$C$ is admissible.
Next we want to show that $(A,C)$ satisfies the Hautus test using
\cref{examples:lemRealHTImpliesHT}.
\begin{remark}
  \label{example2:PoincareImpliesHTForSmallS}
  Let $s_0 \in (0,\infty)$ and $f\in X_1$. Then Poincar\'e's inequality
  shows that for all
  $s \in (0,s_0)$ we have
  \begin{equation*}
    \|f\|^2 \leq   \|f^\prime \|^2 \leq \frac{s_0}{s} \| f^\prime\|^2.
  \end{equation*}
  Therefore we only need to verify the Hautus test estimate for all $s
  \in [s_0,\infty)$ for some fixed $s_0 \in (0,\infty)$.
\end{remark}
The following lemma is one of the key steps in proving that $(A,C)$
satisfies the Hautus test.
\begin{lemma}
  \label{example2:keyIdea}
  Let $s \in (0,\infty)$ and $\eta: (0,1) \to \RR$ be
  absolutely continuous with
  \begin{equation}
    \label{example2:etaAssumptions}
    0 \leq \eta \leq 1, \quad \| \eta^\prime \|_\infty \leq
    \frac{\sqrt{s}}{2}, \quad  V \eta^2
    \leq \frac{s}{2} \eta^2.
  \end{equation}
  Then for all $f \in X_1$ it is true that
  \begin{equation*}
    \| \eta f \|^2 \leq \frac{8}{s^2}\| (B- s\id ) f \|^2
    +\frac{12}{s} \|  C f \|^2.
  \end{equation*}
\end{lemma}
\begin{proof}
  Let $f\in X_1$.
  Note that $\eta^2 f \in \mathcal{D}$, but not necessarily in $X_1$. By
  definition of $B$ and the fact that
  $\eta$ is real-valued we have
  \begin{equation*}
    s \| \eta f  \|^2 = \mathfrak{b} (f, \eta^2 f) - \langle (B-s \id
    ) f ,\eta^2 f \rangle.
  \end{equation*}
  Now by definition it is true that
  \begin{equation*}
    \mathfrak{b} (f, \eta^2 f)  = \int_0^1 f^\prime
    \overline{{(\eta^2 f )}^\prime} \diff
    \lambda + \int_0^1 V f \overline{\eta^2 f} \diff \lambda.
  \end{equation*}
  From the product rule we get
  \begin{equation*}
    {(\eta^2 f )}^\prime = {(\eta^2 )}^\prime f + \eta^2 f^\prime  =
    2 \eta \eta^\prime f
    + \eta^2 f^\prime.
  \end{equation*}
  Using~\eqref{example2:etaAssumptions} and $\eta^2\leq \eta$ this implies
  \begin{equation}
    \label{example2:etaEstim1}
    s \| \eta f\|^2 \leq \| (B- s \id) f\|  \| \eta f\| + \sqrt{s}
    \| \eta f\| \|f^\prime\| +
    \frac{s}{2} \|\eta f\|^2 + \|f^\prime\|^2.
  \end{equation}
  Using the Peter--Paul inequality we get
  \begin{equation*}
    \| (B- s \id) f\|  \| \eta f\| \leq \frac{2}{s} \| (B- s \id)
    f\|^2 + \frac{s}{8} \|\eta f\|^2
  \end{equation*}
  and
  \begin{equation*}
    \sqrt{s} \| \eta f\| \|f^\prime\|  \leq 2 \| f^\prime\|^2 +
    \frac{s}{8} \|\eta f\|^2.
  \end{equation*}
  Combining this with~\eqref{example2:etaEstim1} yields
  \begin{equation*}
    \frac{s}{4} \| \eta f \|^2 \leq \frac{2}{s} \| (B- s \id)
    f\|^2 + 3 \|f^\prime\|^2. \qedhere
  \end{equation*}
\end{proof}
In view of \cref{example2:PoincareImpliesHTForSmallS} and
\cref{example2:keyIdea} it is
sufficient, for proving that $(A,C)$ satisfies the Hautus test, to find $s_0 \in
(0,\infty)$ and constants $c_1, c_2 \in (0,\infty)$ such that for all
$s \in [s_0, \infty)$ there is a function $\eta$ satisfying the assumptions of
\cref{example2:keyIdea} and such that for all $f \in X_1$ we have
\begin{equation*}
  \| f\|^2 \leq c_1 \| \eta f\|^2 +  c_2 s^{-1}\| f^\prime \|^2.
\end{equation*}
The main technical tool to show the existence of such functions is
the following lemma.
\begin{lemma}
  \label{example2:ObsLemma}
  Let $f \in H_0^1 ((0,1))$ and let $I,L,K \subset (0,1)$ be nonempty open
  intervals with $I
  \subset K$ and $L \subset K$. Then we have
  \begin{equation*}
    \int_I |f|^2 \diff \lambda \leq 2
    \frac{\lambda(I)}{\lambda(L)} \int_L |f|^2 \diff \lambda +
    2\lambda (I) \lambda (K) \int_K |f^\prime |^2 \diff \lambda.
  \end{equation*}
\end{lemma}
\begin{proof}
  Let $x \in I$ and $y \in L$. Then the absolute continuity of $f$ implies that
  \begin{equation*}
    | f(x) - f(y) |\leq  \int_{K} | f^\prime|  \diff \lambda
  \end{equation*}
  and therefore it is true that
  \begin{equation*}
    | f(x)|^2 \leq 2 ( |f(x) - f(y)|^2 + |f(y)|^2) \leq
    2\bigg(|f(y)|^2 +\lambda (K)
    \int_{K} | f^\prime|^2  \diff \lambda \bigg).
  \end{equation*}
  By integrating over $L$ and dividing by $\lambda(L)$ we obtain
  \begin{equation*}
    |f(x)|^2 \leq 2 \frac{1}{\lambda(L)} \int_L |f|^2 \diff \lambda +
    2 \lambda (K) \int_K |f^\prime |^2 \diff \lambda.
  \end{equation*}
  Integration over $I$ yields
  \begin{equation*}
    \int_I |f|^2 \diff \lambda \leq 2
    \frac{\lambda(I)}{\lambda(L)} \int_L |f|^2 \diff \lambda +
    2\lambda (I) \lambda (K) \int_K |f^\prime |^2 \diff \lambda. \qedhere
  \end{equation*}
\end{proof}
\begin{proposition}\label{example2:etasProperties}
  For $s \in (0,\infty)$ let $r_s \coloneq 4 s^{-1/2}$,
  \begin{equation*}
    J_s \coloneq \{ j \in \NN_+ \colon H_j > s/2 , \ x_j \geq 3 r_s\},
  \end{equation*}
  \begin{equation*}
    Z_s \coloneq \big[0 , 3 r_s \big] \cup \{ x_j \colon j \in J_s  \},
  \end{equation*}
  and $\eta_s : (0,1) \to [0,1]$,
  \begin{equation*}
    \eta_s (x) \coloneq \min \bigg( 1, \max \bigg( 0,
        \frac{2}{r_s}\dist(x, Z_s)-1
    \bigg) \bigg).
  \end{equation*}

  Let $s \in [256, \infty)$. Then $\eta_s$ is
  absolutely continuous and
  \begin{equation*}
    0 \leq \eta_s \leq 1, \quad \| \eta_s^\prime \|_\infty \leq
    \frac{\sqrt{s}}{2}, \quad  V \eta_s^2
    \leq \frac{s}{2} \eta_s^2.
  \end{equation*}
  Furthermore for all $f \in X_1$ we have
  \begin{equation*}
    \|f\|^2 \leq 5 \|\eta_s f\|^2 + \frac{320}{s} \|f^\prime \|^2.
  \end{equation*}
\end{proposition}
\begin{proof}
  Since $2^8= 256\leq s$, we have $4  r_s \leq 16 \cdot 2^{-4} =  1$.
  Clearly $\eta_s$ is absolutely continuous (even piecewise linear),
  satisfies $0 \leq \eta_s
  \leq 1$ and
  \begin{equation*}
    \| \eta_s^\prime \|_\infty \leq \frac{2}{r_s} = \frac{\sqrt{s}}{2}.
  \end{equation*}

  Let $j \in \NN_+$ such that $H_j>s/2$.
  Then it is true that
  \begin{equation*}
    d_j^2 H_j = 2^{(\ell -2k)j} \leq 1,
  \end{equation*}
  by the assumption that $\ell \leq 2k$. Therefore we have
  \begin{equation*}
    d_j \leq \frac{1}{\sqrt{H_j}} \leq \sqrt{ 2/s} \leq \frac{r_s}{2}.
  \end{equation*}
  Note that for all $x \in (0,1)$, $\eta_s (x) = 0$ whenever $|x-x_j|
  \leq r_s/2$.
  This shows that $\eta_s=0$ on $[x_j - d_j, x_j+d_j]$.
  Therefore $ V \eta_s^2 \leq \frac{s}{2} \eta_s^2 $.

  For $j \in J_s$ let
  \begin{equation*}
    I_j \coloneq (x_j-r_s , x_j+r_s), \quad L_j \coloneq (x_j+r_s,
    x_j+ 2 r_s) ,
    \quad K_j \coloneq (x_j -r_s , x_j +2r_s).
  \end{equation*}
  Let $j \in J_s$. Then $x_j + 2 r_s \leq  x_j +\frac{2}{3} x_j \leq
  \frac{5}{6} \leq 1$ and $x_j -r_s \geq 3 r_s -r_s = 2 r_s \geq 0$.
  That is $K_j \subset (0,1)$.
  Let $j_1, j_2 \in J_s$ with $j_1 <j_2$. Then, using $3 r_s \leq
  x_{j_2}$ we get
  \begin{equation*}
    x_{j_2} + 2r_s  = x_{j_2} + 3r_s -r_s \leq  2 x_{j_2} - r_s \leq
    x_{j_1} -r_s.
  \end{equation*}
  This shows that $K_{j_1} \cap K_{j_2} = \varnothing$.
  Note that $\eta_s =1$ outside of $\bigcup_{j \in J_s} I_j \cup [0,4r_s]$.
  The situation is visualised in \cref{example2:etaPlot}.
  \begin{figure}[htbp]
    \centering
    \includegraphics[width=0.7\textwidth]{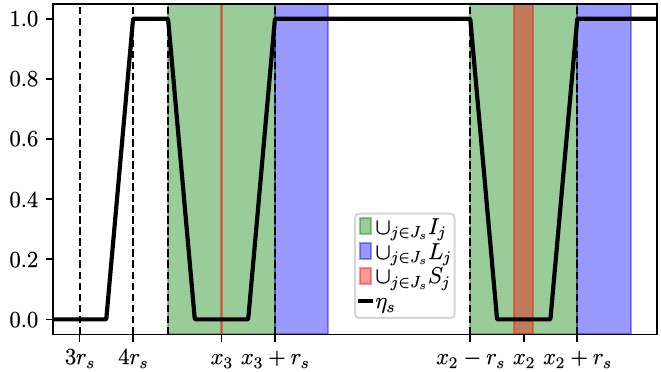}
    \caption{Visualization of the situation in the proof of
      \cref{example2:etasProperties}
      for $k=4$, $\ell = 8$, and $s=2^{15}$: In this case, we have
      $J_s = \{2,3\}$. Here
      $S_j \coloneq [x_j
    -d_j,x_j+d_j]$.}\label{example2:etaPlot}
  \end{figure}
  Let $f \in X_1$. Using Poincar\' e's inequality and
  \cref{example2:ObsLemma} we get
  \begin{equation*}
    \begin{split}
      \|f\|^2 &\leq \|\eta_s f\|^2 + \int_0^{4r_s} |f|^2 \diff \lambda +
      \sum_{j \in J_s} \int_{I_j} |f|^2 \diff \lambda\\
      &\leq \|\eta_s f\|^2 + 8 r_s^2 \|f^\prime\|^2 + \sum_{j \in
      J_s} \bigg(  2
        \frac{\lambda(I_j)}{\lambda(L_j)} \int_{L_j} |f|^2 \diff \lambda +
        2\lambda (I_j) \lambda (K_j) \int_{K_j} |f^\prime |^2 \diff
      \lambda\bigg) \\
      &\leq \|\eta_s f\|^2 + 8 r_s^2 \|f^\prime\|^2 + 4 \| \eta_s f\|^2 + 12
      r_s^2 \| f^\prime\|^2 \\
      &= 5 \| \eta_s f\|^2 + \frac{320}{s} \|f^\prime \|^2. \qedhere
    \end{split}
  \end{equation*}
\end{proof}
\begin{proposition}
  \label{example2:HTsatisfied}
  The system $(A,C)$ satisfies the Hautus test and therefore, as a
  consequence of
  \cref{introduction:mainThm}, $(A,C)$ is exactly observable.
\end{proposition}
\begin{proof}
  Combine \cref{example2:etasProperties} with
  \cref{example2:PoincareImpliesHTForSmallS}, \cref{example2:keyIdea} and
  \cref{examples:lemRealHTImpliesHT}.
\end{proof}
Finally we show that the assumptions of
\cref{examples:RealPartEOSufficiencyEstimate} are
not satisfied.
\begin{proposition}
  \label{examples2:No-ReAsqEstimate}
  There exists no $c\in (0,\infty)$ such that for all $f \in X_1$ it
  is true that
  \begin{equation*}
    \| {(-\re A)}^{1/2} f \| \leq c\| C f \|.
  \end{equation*}
\end{proposition}
\begin{proof}
  Let $\varphi \in C_c^\infty (\RR)$ such that $\supp \varphi
  \subset [3/4, 5/4]$ and $\varphi =1$ on $[15/16, 17/16]$.
  For $j \in \NN_+$ define $f_j : (0,1) \to \CC$ by $f_j(x) \coloneq
  \varphi (x/x_j)$.
  Let $j \in \NN_+$, $j\geq 2$.
  Then $\supp f_j \subset [\frac{3}{4} x_j,\frac{5}{4} x_j] \subset
  (0,1)$ and so $f_j \in X_1$.
  Furthermore we have
  \begin{equation*}
    d_j/x_j = 2^{-j k } 2^j = 2^{-j(k-1)} \leq 2^{-2 j} \leq \frac{1}{4}
  \end{equation*}
  and so
  \begin{equation*}
    x_j + d_j = x_j (1+ d_j/x_j) = x_j (1+2^{-j(k-1)}) \leq x_j (1+
    2^{-4}) =  \frac{17}{16} x_j.
  \end{equation*}
  Similarly we get $x_j -d_j \geq \frac{15}{16} x_j$.
  This shows that $f_j=1$ on $[x_j-d_j,x_j+d_j]$.
  Therefore it is true that
  \begin{equation*}
    \int_0^1 V |f_j|^2 \diff \lambda \geq  H_j d_j \int_{-1}^1 \psi \diff
    \lambda  = 2^{(\ell -k) j} \int_{-1}^1 \psi \diff \lambda.
  \end{equation*}
  The substitution rule yields
  \begin{equation*}
    \| C f_j\|^2  = x_j^{-1}  \| \varphi^\prime \|_{L^2(\RR)}^2 =
    2^{j} \| \varphi^\prime
    \|_{L^2(\RR)}^2.
  \end{equation*}
  This shows (using $\ell \geq k+2$)
  \begin{equation*}
    \begin{split}
      \frac{\langle -A f_j , f_j\rangle}{\| Cf_j\|^2} &=  \frac{
      \mathfrak{b}(f_j,f_j)}{\| C f_j\|^2}  \\
      &\geq 1 + 2^{(\ell -k -1)j}\| \varphi^\prime \|^{-2}_{L^2(\RR)}
      \int_{-1}^1 \psi \diff \lambda \\
      &\geq
      1+ 2^{j} \| \varphi^\prime \|^{-2}_{L^2(\RR)}
      \int_{-1}^1 \psi \diff \lambda   \to \infty , \quad  (j \to
      \infty). \qedhere
    \end{split}
  \end{equation*}
\end{proof}
\begin{remark}\label{example2:generalizationRemark}
  The strategy we used to show that $(A,C)$ satisfies the Hautus test
  may carry over to
  more general diffusion equations with potentials.
  Let $n \in \NN_+$, let $\Omega \subset\mathbb{R}^n$ be a bounded
  domain that has
  sufficiently regular boundary, and let
  $V: \Omega \to
  [0,\infty)$ be measurable and locally bounded.
  Let $X \coloneq L^2(\Omega)$, $Y \coloneq L^2(\Omega, \CC^n)$ and
  $C f \coloneq
  \nabla f$.
  Then the definition of $A$, the relevant properties of $A$,
  \cref{example2:PoincareImpliesHTForSmallS}, and
  \cref{example2:keyIdea}
  carry over with the obvious modifications.
  Therefore it is sufficient, for proving that $(A,C)$ satisfies the
  Hautus test, to find $s_0 \in
  (0,\infty)$ and constants $c_1, c_2 \in (0,\infty)$ such that for all
  $s \in [s_0, \infty)$ there is a weakly differentiable function
  $\eta: \Omega \to \RR$ with
  \begin{equation*}
    0 \leq \eta \leq 1, \quad \| \nabla \eta \|_\infty \leq
    \frac{\sqrt{s}}{2}, \quad  V \eta^2
    \leq \frac{s}{2} \eta^2,
  \end{equation*}
  and such that for all $f \in X_1$ we have
  \begin{equation*}
    \| f\|^2 \leq c_1 \| \eta f\|^2 +  c_2 s^{-1}\| \nabla f \|^2.
  \end{equation*}
  Importantly the properties $\eta$ needs to satisfy only depend on
  $s$, $V$ and, in view of
  Poincar\' e's inequality, the geometry of $\Omega$.
\end{remark}

\section*{Acknowledgments}
ChatGPT was used to proofread the paper, to find
relevant literature, and to provide some initial proof ideas.
Nevertheless, the final proofs in this paper are developed by the
authors, using the
ChatGPT-generated proof ideas as a starting point
when applicable.
The authors assume responsibility for all content.
\printbibliography%

@book{TucsnakWeiss,
  author = {Tucsnak, Marius and Weiss, George},
  title = {Observation and control for operator semigroups},
  fseries = {Birkh{\"a}user Advanced Texts. Basler Lehrb{\"u}cher},
  series = {Birkh{\"a}user Adv. Texts, Basler Lehrb{\"u}ch.},
  issn = {1019-6242},
  isbn = {978-3-7643-8993-2; 978-3-7643-8994-9},
  year = {2009},
  publisher = {Birkh{\"a}user},
  language = {English},
  zbMATH = {5359001},
  Zbl = {1188.93002}
}

@article{RussellWeissGeneralCondExactObs,
  title = {A general necessary condition for exact observability},
  volume = {32},
  issn = {0363-0129},
  doi = {10.1137/S036301299119795X},
  pages = {1--23},
  number = {1},
  journaltitle = {{SIAM} Journal on Control and Optimization},
  journal = {{SIAM} J. Control Optim.},
  author = {Russell, David L. and Weiss, George},
  year = {1994},
}

@article{jacobZwartExactObsFinDimOutput,
  title = {Exact observability of diagonal systems with a
  finite-dimensional output operator},
  volume = {43},
  issn = {0167-6911},
  doi = {10.1016/S0167-6911(00)00117-1},
  pages = {101--109},
  number = {2},
  journaltitle = {Systems \& Control Letters},
  journal = {Syst. Control Lett.},
  author = {Jacob, Birgit and Zwart, Hans},
  year = {2001},
}

@article{jacobZwartHautusForExpStable,
  title = {On the {H}autus test for exponentially stable \({C_0}\)-groups},
  volume = {48},
  issn = {0363-0129},
  doi = {10.1137/080724733},
  pages = {1275--1288},
  number = {3},
  journaltitle = {{SIAM} Journal on Control and Optimization},
  journal = {{SIAM} J. Control Optim.},
  author = {Jacob, Birgit and Zwart, Hans},
  year = {2009},
}

@article{jacob_counterexamples_2004,
	title = {Counterexamples concerning observation operators for \({C_0}\)-semigroups},
	volume = {43},
	issn = {0363-0129},
	doi = {10.1137/S0363012903423235},
	pages = {137--153},
	number = {1},
	journaltitle = {{SIAM} Journal on Control and Optimization},
	journal = {{SIAM} J. Control Optim.},
	author = {Jacob, Birgit and Zwart, Hans},
	year = {2004},
}

@article{miller_controllability_2005,
	title = {Controllability cost of conservative systems: resolvent condition and transmutation},
	volume = {218},
	issn = {0022-1236},
	doi = {10.1016/j.jfa.2004.02.001},
	pages = {425--444},
	number = {2},
	journaltitle = {Journal of Functional Analysis},
	journal = {J. Funct. Anal.},
	author = {Miller, Luc},
	year = {2005},
}

@book{rudin_functional_1991,
	edition = {2nd ed},
	title = {Functional analysis},
	isbn = {0-07-054236-8},
	publisher = {{McGraw}-Hill},
	author = {Rudin, Walter},
	year = {1991},
}

@incollection{jacob_zwart_blondel_unsolved_2004,
booktitle = {Unsolved problems in mathematical systems and control theory},
	isbn = {0-691-11748-9},
	language = {English},
	publisher = {Princeton University Press},
	editor = {Blondel, Vincent D. and Megretski, Alexandre},
	year = {2004},
	author = {Jacob, Birgit and Zwart, Hans},
  title   = {A {H}autus test for infinite-dimensional systems},
pages = {251--255},
}

@book{hytonen_analysis_2016,
	series = {Ergeb. {Math}. {Grenzgeb}., 3. {Folge}},
	title = {Analysis in {Banach} spaces. {Volume} {I}. {Martingales} and {Littlewood}-{Paley} theory},
	volume = {63},
	isbn = {978-3-319-48519-5 978-3-319-48520-1},
	issn = {0071-1136},
	doi = {10.1007/978-3-319-48520-1},
	language = {English},
	publisher = {Springer},
	author = {Hytönen, Tuomas and van Neerven, Jan and Veraar, Mark and Weis, Lutz},
	year = {2016},
}

@misc{daun2026hautustestexactobservability,
      title={On the Hautus test for exact observability of normal semigroups}, 
      author={Jannik Daun and Birgit Jacob},
      year={2026},
      eprint={2608.09407},
      archivePrefix={arXiv},
      primaryClass={math.OC},
      url={https://arxiv.org/abs/2608.09407}, 
}

@book{cohn_measure_2013,
	edition = {2nd revised ed.},
	series = {Birkhäuser {Adv}. {Texts}, {Basler} {Lehrbüch}.},
	title = {Measure theory},
	isbn = {978-1-4614-6955-1 978-1-4614-6956-8},
	issn = {1019-6242},
	doi = {10.1007/978-1-4614-6956-8},
	language = {English},
	publisher = {Birkhäuser},
	author = {Cohn, Donald L.},
	year = {2013},
}

@book{kato_perturbation_1995,
	edition = {Reprint of the corr. print. of the 2nd ed. 1980},
	title = {Perturbation theory for linear operators},
	isbn = {3-540-58661-X},
	issn = {1431-0821},
	doi = {10.1007/978-3-642-66282-9},
	series = {Class. Math.},
	publisher = {Berlin: Springer-Verlag},
	author = {Kato, Tosio},
	date = {1995},
}
\end{document}